\documentclass[12pt,a4paper]{article}

\usepackage{theorem,enumerate}
\usepackage{amsmath,latexsym,amssymb,amsfonts}
\usepackage{eucal}
\usepackage{eufrak}

\newcounter{Scounter}
\theorembodyfont{\normalfont\slshape}
\newtheorem{thm}{Theorem}

\newtheorem{prop}[thm]{Proposition}
\newtheorem{lem}[thm]{Lemma}

\newtheorem{Q}[thm]{Question} 

\newtheorem{con}[thm]{Conjecture}

\newcommand{\proof}{\medbreak\noindent\textit{Proof.}\quad}

\newcommand{\qed}{{$\quad\square$\vs{3.6}}}

\newcommand{\vs}[1]{\vspace*{#1 mm}}

\numberwithin{equation}{section}

\def\L{{ \mathcal{L}}}

\def\P{{ \mathcal{P}}}

\def\Q{{ \mathcal{Q}}}

\bfseries\normalfont

\date{}

\begin{document}
\title{A New Approach Towards a Conjecture on Intersecting Three Longest Paths}
\author{Shinya Fujita$^1$,  Michitaka Furuya$^2$, Reza Naserasr$^3$, Kenta Ozeki$^4$}

\maketitle

\noindent\footnotetext[1]{International College of Arts and Sciences, Yokohama City University,  22-2 Seto, Kanazawa-ku, Yokohama 236-0027, Japan; shinya.fujita.ph.d@gmail.com}
\noindent\footnotetext[2]{Department of Mathematical Information Science, Tokyo University of Science, 1-3 Kagurazaka, Shinjuku-ku, Tokyo 162-8601, Japan; michitaka.furuya@gmail.com}
\noindent\footnotetext[3]{CNRS, LRI, UMR8623, Univ. Paris-Sud 11, F-91405 Orsay Cedex, France; reza@lri.fr}
\noindent\footnotetext[4]{National Institute of Informatics, 2-1-2 Hitotsubashi, Chiyoda-Ku, Tokyo 101-8430, Japan and JST, ERATO, Kawarabayashi Large Graph Project, Japan; ozeki@nii.ac.jp}

\begin{abstract}
In 1966, T. Gallai asked whether  every connected graph has a vertex that appears in all longest paths. Since then this question has attracted much attention and many work has been done in this topic. 
One important open question in this area is to ask whether any three longest paths contains a common vertex in a connected graph. 
It was conjectured that the answer to this question is positive. 
In this paper, we propose a new approach in view of distances among longest paths in a connected graph,
 and give a substantial progress towards the conjecture along the idea.   
\end{abstract}
\section{Introduction}\label{sec1}

In \cite{ref5} Gallai asked whether every connected graph has a vertex that appears in all longest paths. This question has attracted much attention and many work has been done around this area of study. 
The answer to this question is false as stated; actually several counterexamples were given in \cite{ref16, ref18, ref21}.  
A graph $G$ is \textit{hypotraceable} if $G$ has no Hamiltonian path but every vertex-deleted subgraph $G-v$ has. 
Note that hypotraceable graphs constitute a large class of counterexamples. 
 Thomassen \cite{ref12} showed that there exist infinitely many planar hypotraceable graphs, meaning that there exist 
 infinitely many counterexamples towards the question. 
 
 Yet there are classes of graphs for which the answer to Gallai's question is positive. 
 To see this, note that, in a tree, all longest paths must contain its center(s). 
 Klav\u{z}ar and Petkov\u{s}ek \cite{ref9} showed that the answer is also positive for split graphs, cacti, and some other classes of graphs. 
 Balister et al. \cite{ref2} obtained a similar result for the class of circular arc graphs. 
 
Regarding Gallai's question, what happens if we consider the intersection of a smaller number of longest paths? 
While we can easily check that every two longest paths share a vertex, it is not known whether every three longest paths also share a vertex. 
In \cite{ref7} it appears as a conjecture, which has originally been asked by Zamfirescu since the 1980s (see \cite{ref22}). 

\begin{con}
\label{Z}
For every connected graph, any three of its longest paths have a common vertex. 
\end{con}

So far, very little progress has been made on this conjecture. 
Axenovich \cite{ref1} proved that Conjecture~\ref{Z} is true for connected outerplanar graphs, and de Rezende et al. \cite{ref0} proved that Conjecture~\ref{Z} is true for connected graphs in which all nontrivial blocks are Hamiltoninan.  

In this paper, we introduce a new graph parameter in view of distances among longest paths in a connected graph. 
To state this, we give some basic definitions. 
For a graph $G$, let $P$ be a path in $G$, and let $x$ and $y$ be the end-vertices of $P$.
Note that $|V(P)|=1$ if and only if $x=y$.
For $X,Y\subseteq V(G)$, $P$ is called an {\it $X$-$Y$ path} if $V(P)\cap X=\{x\}$ and $V(P)\cap Y=\{y\}$.
Let $u,v\in V(P)$.
We let $uPv$ denote the $\{u\}$-$\{v\}$ path on $P$.
Furthermore, we let $\check{u}Pv=uPv-u$, $uP\check{v}=uPv-v$ and $\check{u}P\check{v}=uPv-\{u,v\}$.

Let $G$ be a connected graph. 
Let $l(G)$ be the length of any longest path in $G$, and let $\L(G)$ be the set of longest paths of $G$; thus $\L(G)=\{P\mid P$ is a path in $G$ with $|V(P)|=l(G)+1\}$. 
For $x,y\in V(G)$ let $d_G(x,y)$ be the distance between $x$ and $y$ in $G$ (i.e., the length of a shortest path joining $x$ and $y$ in $G$). 
Also for a vertex $x\in V (G)$ and a subset $U\subseteq V(G)$, let $d_G(x,U) = \min\{d_G(x,y) \mid y\in U\}$.
For $\P\subseteq \L(G)$, let $f(G,\P)=\min \{\sum_{P\in \P}d_{G}(v,V(P))\mid v\in V(G)\}$. 

Using this graph parameter, we can formulate Conjecture~\ref{Z} as follows.

\begin{con}
\label{Gallai}
Let $G$ be a connected graph, and let $\P$ be a subset of $\L(G)$ with $|\P|=3$.
Then $f(G,\P)=0$.
\end{con}

As mentioned before, it is easy to check that any two longest paths of a connected graph have a common vertex. We now give the proof in this context. 

\begin{prop}
\label{prop1}
Let $G$ be a connected graph, and let $\P$ be a subset of $\L(G)$ with $|\P|=2$.
Then $f(G,\P)=0$.
\end{prop}
\proof
Write $\P=\{P_{1},P_{2}\}$, and for each $i\in \{1,2\}$, let $u_{i}$ and $v_{i}$ be the end-vertices of $P_{i}$.
Since $G$ is connected, $G$ has a $V(P_{1})$-$V(P_{2})$ path $Q$.
Note that $V(P_{1})\cap V(P_{2})\not=\emptyset $ if and only if $|V(Q)|=1$.
For each $i\in \{1,2\}$, write $V(P_{i})\cap V(Q)=\{w_{i}\}$.
We may assume that $|V(u_{i}P_{i}w_{i})|\geq |V(v_{i}P_{i}w_{i})|$ for each $i\in \{1,2\}$.
Then the length of the path $u_{1}P_{1}w_{1}Qw_{2}P_{2}u_{2}$ in $G$ is $(|V(u_{1}P_{1}w_{1})|-1)+(|V(Q)|-1)+(|V(u_{2}P_{2}w_{2})|-1)$.
On the other hand, for each $i\in \{1,2\}$, $|V(u_{i}P_{i}w_{i})|-1\geq ((|V(u_{i}P_{i}w_{i})|-1)+(|V(v_{i}P_{i}w_{i})|-1))/2=(|V(P_{i})|-1)/2=l(G)/2$.
Consequently,
\begin{eqnarray}
\begin{split}
(|V(u_{1}P_{1}w_{1})|-1)&+(|V(u_{2}P_{2}w_{2})|-1)\nonumber \\
&\geq \frac{l(G)}{2}+\frac{l(G)}{2}\nonumber \\
&= l(G)\nonumber \\
&\geq (|V(u_{1}P_{1}w_{1})|-1)+(|V(Q)|-1)+(|V(u_{2}P_{2}w_{2})|-1).\nonumber
\end{split}
\end{eqnarray}
This leads to $|V(Q)|=1$, and hence $V(P_{1})\cap V(P_{2})\not=\emptyset $.
\qed

In this paper, we give an upper bound of $f(G,\P)$ with $|\P|=3$, which is linear in terms of $|V(G)|$.

\begin{thm}
\label{thm1}
Let $G$ be a connected graph of order $n$, and let $\P$ be a subset of $\L(G)$ with $|\P|=3$.
Then $f(G,\P)\leq (n+6)/13$.
\end{thm}

After proving this bound in Section 2, in the follow-up section we show that to 
prove the conjecture it would be enough to improve our linear bound to any 
nondecreasing sublinear bound.
Namely, we propose an equivalent conjecture towards Conjecture~\ref{Z} in terms of the function $f(G,\P)$.

\section{Proof of Theorem~\ref{thm1}}\label{sec2}

We start with some lemmas.

For a set $\P$ of graphs and $P\in \P$, set $X_{\P}(P)=V(P)-(\bigcup _{P'\in \P-\{P\}}V(P'))$.

\begin{lem}
\label{lem2.1}
Let $G$ be a connected graph of order $n$, and let $\P\subseteq \L(G)$ with $|\P|=3$.
If $f(G,\P)>0$, then $n\geq (3l(G)+\sum _{P\in \P}|X_{\P}(P)|+3)/2$.
\end{lem}
\proof
Write $\P=\{P_{1},P_{2},P_{3}\}$.
Since $\bigcap _{1\leq i\leq 3}V(P_{i})=\emptyset $,
\begin{align}
n\geq |\bigcup _{1\leq i\leq 3}V(P)|=\sum _{1\leq i\leq 3}|X_{\P}(P_{i})|+\sum _{1\leq i<j\leq 3}|V(P_{i})\cap V(P_{j})|.\label{eq2.1}
\end{align}
Since $l(G)+1=|V(P_{i})|=|X_{\P}(P_{i})|+\sum _{j\not=i}|V(P_{i})\cap V(P_{j})|$ for each $1\leq i\leq 3$,
\begin{eqnarray}
3l(G)+3 &=& \sum _{1\leq i\leq 3}|X_{\P}(P_{i})|+\sum _{1\leq i\leq 3}(\sum _{j\not=i}|V(P_{i})\cap V(P_{j})|)\nonumber \\
&=& \sum _{1\leq i\leq 3}|X_{\P}(P_{i})|+2\sum _{1\leq i<j\leq 3}|V(P_{i})\cap V(P_{j})|.\label{eq2.2}
\end{eqnarray}
By (\ref{eq2.1}) and (\ref{eq2.2}),
\begin{eqnarray}
n &\geq & \sum _{1\leq i\leq 3}|X_{\P}(P_{i})|+\sum _{1\leq i<j\leq 3}|V(P_{i})\cap V(P_{j})|\nonumber \\
&=& \sum _{1\leq i\leq 3}|X_{\P}(P_{i})|+(3l(G)+3-\sum _{1\leq i\leq 3}|X_{\P}(P_{i})|)/2\nonumber \\
&=& (3l(G)+3+\sum _{1\leq i\leq 3}|X_{\P}(P_{i})|)/2.\nonumber
\end{eqnarray}
Thus we get the desired conclusion.
\qed

For a set $\P$ of three paths and $P\in \P$, let $t_{\P}(P)$ be the number of $V(P_{1})$-$V(P_{2})$ paths on $P$, where $\P-\{P\}=\{P_{1},P_{2}\}$.
If $\P$ consists of three longest paths of a connected graph, then $t_{\P}(P)\geq 1$ for every $P\in \P$ by Proposition~\ref{prop1}.

\begin{lem}
\label{lem2.2}
Let $G$ be a connected graph, and let $\P\subseteq \L(G)$ with $|\P|=3$.
Then $|X_{\P}(P)|\geq t_{\P}(P)(f(G,\P)-1)$ for each $P\in \P$.
\end{lem}
\proof
We may assume that $f(G,\P)\geq 1$.
Write $\P-\{P\}=\{P_{1},P_{2}\}$, and let $\Q$ be the set of $V(P_{1})$-$V(P_{2})$ paths on $P$.
Note that every path in $\Q$ has order at least two and $|\Q|=t_{\P}(P)$.
Let $Q\in \Q$, and let $u$ and $v$ be the end-vertices of $Q$ with $u\in V(P_{1})$ and $v\in V(P_{2})$.
Then $V(Q)\cap X_{\P}(P)=V(Q)-\{u,v\}$.
Since $u\in V(P)\cap V(P_{1})$, $f(G,\P)\leq \sum_{P'\in \P}d_{G}(u,V(P'))=d_{G}(u,V(P_{2}))\leq d_{G}(u,v)\leq |V(Q)|-1$.
Hence $|V(Q)\cap X_{\P}(P)|=|V(Q)|-2\geq f(G,\P)-1$.
Since $Q$ is arbitrary,
\begin{align}
\sum _{Q\in \Q}|V(Q)\cap X_{\P}(P)|\geq t_{\P}(P)(f(G,\P)-1).\label{eq2.3}
\end{align}
Clearly, each vertex in $X_{\P}(P)$ belongs to at most one path in $\Q$.
This together with (\ref{eq2.3}) implies that $|X_{\P}(P)|\geq |\bigcup _{Q\in \Q}(V(Q)\cap X_{\P}(P))|=\sum _{Q\in \Q}|V(Q)\cap X_{\P}(P)|\geq t_{\P}(P)(f(G,\P)-1)$.
\qed

\begin{lem}
\label{lem2.3}
Let $G$ be a connected graph, and let $\P\subseteq \L(G)$ with $|\P|=3$.
If there exists a path $P\in \P$ with $t_{\P}(P)=1$, then $f(G,\P)=0$.
\end{lem}
\proof
Suppose that $f(G,\P)>0$.
Let $u$ and $v$ be the end-vertices of $P$.
Write $\P-\{P\}=\{P_{1},P_{2}\}$, and for each $i\in \{1,2\}$, let $w_{i}$ be the vertex which is contained in $P_{i}$ and the unique $V(P_{1})$-$V(P_{2})$ path on $P$ (see Figure~\ref{fig1}).
We may assume that $|V(uPw_{1})|\leq |V(uPw_{2})|$.
Since $f(G,\P)>0$, $w_{1}\not=w_{2}$, and hence $|V(w_{1}Pv)|>|V(w_{2}Pv)|$.
Furthermore, we may assume that $|V(uPw_{1})|\leq |V(vPw_{2})|$.
Since $l(G)=|V(uPw_{1})|+|V(w_{1}Pv)|-2$,
\begin{eqnarray}
|V(w_{1}Pv)| &>& \frac{|V(w_{1}Pv)|}{2}+\frac{|V(w_{2}Pv)|}{2}\nonumber \\
&=& \frac{l(G)-|V(uPw_{1})|+2}{2}+\frac{|V(w_{2}Pv)|}{2}\nonumber \\
&\geq & \frac{l(G)-|V(vPw_{2})|+2}{2}+\frac{|V(w_{2}Pv)|}{2}\nonumber \\
&=& \frac{l(G)+2}{2}.\label{eq2.4}
\end{eqnarray}
Let $u_{1}$ and $v_{1}$ be the end-vertices of $P_{1}$.
We may assume that $|V(u_{1}P_{1}w_{1})|\geq |V(w_{1}P_{1}v_{1})|$.
Since $l(G)=|V(u_{1}P_{1}w_{1})|+|V(w_{1}P_{1}v_{1})|-2$,
\begin{align}
|V(u_{1}P_{1}w_{1})|\geq \frac{|V(u_{1}P_{1}w_{1})|+|V(w_{1}P_{1}v_{1})|}{2}=\frac{l(G)+2}{2}.\label{eq2.4+}
\end{align}
By (\ref{eq2.4}) and (\ref{eq2.4+}), $|V(u_{1}P_{1}w_{1})|+|V(w_{1}Pv)|-2>(l(G)+2)/2+(l(G)+2)/2-2=l(G)$.
By the assumption that $t_{\P}(P)=1$, the path $\check{w}_{1}Pv$ contains no vertex in $V(P_{1})$.
Hence $P^{(1)}_{1}=u_{1}P_{1}w_{1}Pv$ is a path in $G$ with length $|V(u_{1}P_{1}w_{1})|+|V(w_{1}Pv)|-2>l(G)$, which is a contradiction.
\qed

\begin{figure}
\begin{center}
\unitlength 0.1in
\begin{picture}( 23.6500, 14.7200)(  6.8300,-18.2700)
%
\special{sh 1.000}%
\special{ia 774 1092 36 36  0.0000000 6.2831853}%
\special{pn 8}%
\special{ar 774 1092 36 36  0.0000000 6.2831853}%
%
\special{sh 1.000}%
\special{ia 1054 1092 36 36  0.0000000 6.2831853}%
\special{pn 8}%
\special{ar 1054 1092 36 36  0.0000000 6.2831853}%
%
\special{sh 1.000}%
\special{ia 1614 1092 36 36  0.0000000 6.2831853}%
\special{pn 8}%
\special{ar 1614 1092 36 36  0.0000000 6.2831853}%
%
\special{sh 1.000}%
\special{ia 2174 1092 36 36  0.0000000 6.2831853}%
\special{pn 8}%
\special{ar 2174 1092 36 36  0.0000000 6.2831853}%
%
\special{sh 1.000}%
\special{ia 2734 1092 36 36  0.0000000 6.2831853}%
\special{pn 8}%
\special{ar 2734 1092 36 36  0.0000000 6.2831853}%
%
\special{sh 1.000}%
\special{ia 3014 1092 36 36  0.0000000 6.2831853}%
\special{pn 8}%
\special{ar 3014 1092 36 36  0.0000000 6.2831853}%
%
\special{pn 8}%
\special{pa 3014 1092}%
\special{pa 2594 1092}%
\special{fp}%
\special{pa 2314 1092}%
\special{pa 2034 1092}%
\special{fp}%
\special{pa 1754 1092}%
\special{pa 1474 1092}%
\special{fp}%
\special{pa 1194 1092}%
\special{pa 774 1092}%
\special{fp}%
%
\special{pn 4}%
\special{sh 1}%
\special{ar 1894 1092 6 6 0  6.28318530717959E+0000}%
\special{sh 1}%
\special{ar 1824 1092 6 6 0  6.28318530717959E+0000}%
\special{sh 1}%
\special{ar 1964 1092 6 6 0  6.28318530717959E+0000}%
\special{sh 1}%
\special{ar 1964 1092 6 6 0  6.28318530717959E+0000}%
%
\special{pn 4}%
\special{sh 1}%
\special{ar 2454 1092 6 6 0  6.28318530717959E+0000}%
\special{sh 1}%
\special{ar 2384 1092 6 6 0  6.28318530717959E+0000}%
\special{sh 1}%
\special{ar 2524 1092 6 6 0  6.28318530717959E+0000}%
\special{sh 1}%
\special{ar 2524 1092 6 6 0  6.28318530717959E+0000}%
%
\special{pn 4}%
\special{sh 1}%
\special{ar 1334 1092 6 6 0  6.28318530717959E+0000}%
\special{sh 1}%
\special{ar 1264 1092 6 6 0  6.28318530717959E+0000}%
\special{sh 1}%
\special{ar 1404 1092 6 6 0  6.28318530717959E+0000}%
\special{sh 1}%
\special{ar 1404 1092 6 6 0  6.28318530717959E+0000}%
\put(7.7300,-12.3200){\makebox(0,0){$u$}}%
\put(30.1300,-12.3200){\makebox(0,0){$v$}}%
\put(9.1300,-4.2000){\makebox(0,0){$u_{1}$}}%
\put(17.5300,-17.9200){\makebox(0,0){$v_{1}$}}%
\put(17.4000,-11.9000){\makebox(0,0){$w_{1}$}}%
\put(22.9800,-11.9700){\makebox(0,0){$w_{2}$}}%
\put(28.7300,-9.9400){\makebox(0,0){$P$}}%
\put(15.9200,-6.7200){\makebox(0,0){$P_{1}$}}%
\put(24.5300,-8.2600){\makebox(0,0){$P_{2}$}}%
%
\special{pn 8}%
\special{pa 914 532}%
\special{pa 920 580}%
\special{pa 926 628}%
\special{pa 934 674}%
\special{pa 940 720}%
\special{pa 948 762}%
\special{pa 956 802}%
\special{pa 964 838}%
\special{pa 974 872}%
\special{pa 986 898}%
\special{pa 996 922}%
\special{pa 1010 940}%
\special{pa 1024 950}%
\special{pa 1040 954}%
\special{pa 1056 952}%
\special{pa 1076 942}%
\special{pa 1096 924}%
\special{pa 1116 902}%
\special{pa 1138 876}%
\special{pa 1162 848}%
\special{pa 1210 788}%
\special{pa 1234 758}%
\special{pa 1258 730}%
\special{pa 1282 708}%
\special{pa 1306 688}%
\special{pa 1328 674}%
\special{pa 1352 668}%
\special{pa 1372 668}%
\special{pa 1394 676}%
\special{pa 1414 688}%
\special{pa 1434 706}%
\special{pa 1452 728}%
\special{pa 1470 754}%
\special{pa 1488 784}%
\special{pa 1506 818}%
\special{pa 1522 854}%
\special{pa 1538 892}%
\special{pa 1554 932}%
\special{pa 1570 972}%
\special{pa 1600 1054}%
\special{pa 1614 1094}%
\special{pa 1642 1170}%
\special{pa 1654 1206}%
\special{pa 1668 1242}%
\special{pa 1682 1276}%
\special{pa 1694 1308}%
\special{pa 1706 1342}%
\special{pa 1718 1374}%
\special{pa 1754 1464}%
\special{pa 1766 1492}%
\special{pa 1778 1520}%
\special{pa 1790 1548}%
\special{pa 1812 1602}%
\special{pa 1834 1654}%
\special{pa 1856 1704}%
\special{pa 1866 1730}%
\special{pa 1876 1754}%
\special{pa 1888 1778}%
\special{pa 1894 1792}%
\special{sp}%
%
\special{sh 1.000}%
\special{ia 914 532 36 36  0.0000000 6.2831853}%
\special{pn 8}%
\special{ar 914 532 36 36  0.0000000 6.2831853}%
%
\special{sh 1.000}%
\special{ia 1894 1792 36 36  0.0000000 6.2831853}%
\special{pn 8}%
\special{ar 1894 1792 36 36  0.0000000 6.2831853}%
%
\special{pn 8}%
\special{pa 2314 812}%
\special{pa 2034 1372}%
\special{fp}%
%
\special{pn 4}%
\special{sh 1}%
\special{ar 1964 1512 6 6 0  6.28318530717959E+0000}%
\special{sh 1}%
\special{ar 1998 1442 6 6 0  6.28318530717959E+0000}%
\special{sh 1}%
\special{ar 1928 1582 6 6 0  6.28318530717959E+0000}%
\special{sh 1}%
\special{ar 1928 1582 6 6 0  6.28318530717959E+0000}%
%
\special{pn 4}%
\special{sh 1}%
\special{ar 2384 672 6 6 0  6.28318530717959E+0000}%
\special{sh 1}%
\special{ar 2418 602 6 6 0  6.28318530717959E+0000}%
\special{sh 1}%
\special{ar 2348 742 6 6 0  6.28318530717959E+0000}%
\special{sh 1}%
\special{ar 2348 742 6 6 0  6.28318530717959E+0000}%
\end{picture}%
\caption{paths in $\P$}
\label{fig1}
\end{center}
\end{figure}

\medbreak\noindent\textit{Proof of Theorem~\ref{thm1}.}\quad
We may assume that $f(G,\P)\geq 1$.
Choose $P\in \P$ so that $t=t_{\P}(P)$ is as small as possible.
Then $t_{\P}(P)\geq 2$ by Lemma~\ref{lem2.3}.
Let $u$ and $v$ be the end-vertices of $P$.
Write $\P-\{P\}=\{P_{1},P_{2}\}$, and let $u_{i}$ and $v_{i}$ be the end-vertices of $P_{i}$ for each $i\in \{1,2\}$.
Let $Q_{1},Q_{2},\cdots ,Q_{t}$ be the $V(P_{1})$-$V(P_{2})$ paths on $P$ which are aligned on $P$ in order of indices with initial point $u$ (i.e. for each $2\leq i\leq t$, the unique $\{u\}$-$V(Q_{i})$ path on $P$ contains $\bigcup _{1\leq j\leq i-1}V(Q_{j})$).
We may assume that the length of the unique $\{u\}$-$V(Q_{1})$ path on $P$ is at least that of the unique $\{v\}$-$V(Q_{t})$ path on $P$.
For each $1\leq i\leq t$ and each $j\in \{1,2\}$, write $V(Q_{i})\cap V(P_{j})=\{w^{(j)}_{i}\}$.
We may assume that $|V(uPw^{(1)}_{1})|\leq |V(uPw^{(2)}_{1})|$.
Let $R$ be a $\{w^{(1)}_{1}\}$-$V(P_{2})$ path on $P_{1}$, and write $V(R)\cap V(P_{2})=\{x\}$.
For each $i\in \{1,2\}$, we may assume that $|V(u_{i}P_{i}w^{(i)}_{1})|\leq |V(u_{i}P_{i}x)|$ (see Figure~\ref{fig2}).

\begin{figure}
\begin{center}
\input{f2.tex}
\caption{paths in $\P$}
\label{fig2}
\end{center}
\end{figure}

Since $w^{(1)}_{1}\in V(P)\cap V(P_{1})$, $f(G,\P)\leq \sum_{P'\in \P}d_{G}(w^{(1)}_{1},V(P'))=d_{G}(w^{(1)}_{1},V(P_{2}))\leq \min \{d_{G}(w^{(1)}_{1},w^{(2)}_{1}),d_{G}(w^{(1)}_{1},x)\}\leq \min \{|V(Q_{1})|-1,|V(R)|-1\}$.
Hence
\begin{align}
|V(Q_{1})|\geq f(G,\P)+1~\mbox{ and }~|V(R)|\geq f(G,\P)+1.\label{eq2.5}
\end{align}

\begin{figure}
\begin{center}
\input{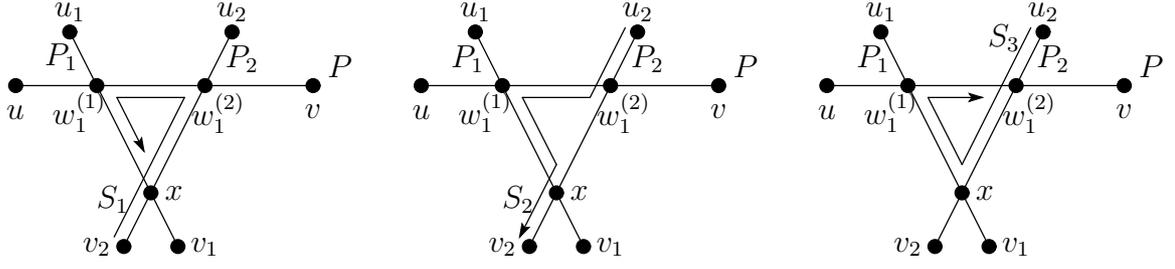}
\caption{path $S_{i}$}
\label{fig3}
\end{center}
\end{figure}

Since $w^{(2)}_{1}Q_{1}\check{w}^{(1)}_{1}$ contains no vertex in $V(P_{1})$, $w^{(2)}_{1}Q_{1}w^{(1)}_{1}Rx$ is a path in $G$.
Furthermore, since $\check{w}^{(2)}_{1}Q_{1}w^{(1)}_{1}P_{1}\check{x}$ contains no vertex in $V(P_{2})$,\\
(i)~$S_{1}=v_{2}P_{2}w^{(2)}_{1}Q_{1}w^{(1)}_{1}R\check{x}$,\\
(ii)~$S_{2}=u_{2}P_{2}w^{(2)}_{1}Q_{1}w^{(1)}_{1}RxP_{2}v_{2}$ and\\
(iii)~$S_{3}=u_{2}P_{2}xRw^{(1)}_{1}Q_{1}\check{w}^{(2)}_{1}$.\\
are paths in $G$ (see Figure~\ref{fig3}).

Since the length of $S_{1}$ is $(|V(v_{2}P_{2}w^{(2)}_{1}|-1)+(|V(Q_{1})|-1)+(|V(w^{(1)}_{1}R\check{x})|-1)$ and $|V(w^{(1)}_{1}R\check{x})|=|V(R)|-1$, we have $(|V(v_{2}P_{2}w^{(2)}_{1})|-1)+(|V(w^{(2)}_{1}P_{2}u_{2})|-1)=|V(P_{2})|-1=l(G)\geq (|V(v_{2}P_{2}w^{(2)}_{1}|-1)+(|V(Q_{1})|-1)+(|V(R)|-2)$.
This together with (\ref{eq2.5}) leads to
\begin{align}
|V(u_{2}P_{2}w^{(2)}_{1})|\geq |V(Q_{1})|+|V(R)|-2\geq 2f(G,\P).\label{eq2.6}
\end{align}
By comparing the length of $P_{2}$ and $S_{2}$ and (\ref{eq2.5}), we have
\begin{align}
|V(w^{(2)}_{1}P_{2}x)|\geq |V(Q_{1})|+|V(R)|-1\geq 2f(G,\P)+1.\label{eq2.7}
\end{align}
By comparing the length of $P_{2}$ and $S_{3}$ and (\ref{eq2.5}), we also have
\begin{align}
|V(xP_{2}v_{2})|\geq |V(Q_{1})|+|V(R)|-2\geq 2f(G,\P).\label{eq2.8}
\end{align}
Therefore
\begin{eqnarray}
l(G) &=& |V(P_{2})|-1\nonumber \\
&=& |V(u_{2}P_{2}w^{(2)}_{1})|+|V(w^{(2)}_{1}P_{2}x)|+|V(xP_{2}v_{2})|-3\nonumber \\
&\geq & 2f(G,\P)+(2f(G,\P)+1)+2f(G,\P)-3\nonumber \\
&=& 6f(G,\P)-2.\label{eq2.9}
\end{eqnarray}

\medskip
\noindent
\textbf{Case 1:} $t_{\P}(P)=2$.

\begin{figure}
\begin{center}
\unitlength 0.1in
\begin{picture}( 25.2000, 13.2200)(  4.8300,-16.9200)
\put(5.7300,-9.6000){\makebox(0,0){$u$}}%
\put(9.2400,-9.6700){\makebox(0,0){$w^{(1)}_{1}$}}%
\put(16.3200,-9.6000){\makebox(0,0){$w^{(2)}_{1}$}}%
\put(16.9300,-4.3500){\makebox(0,0){$u_{2}$}}%
\put(8.5300,-4.3500){\makebox(0,0){$u_{1}$}}%
\put(9.9300,-16.6000){\makebox(0,0){$v_{2}$}}%
\put(15.5300,-16.6000){\makebox(0,0){$v_{1}$}}%
\put(13.9200,-13.8000){\makebox(0,0){$x$}}%
\put(17.1400,-6.7700){\makebox(0,0){$P_{2}$}}%
\put(8.3200,-6.7700){\makebox(0,0){$P_{1}$}}%
\put(29.5300,-9.5300){\makebox(0,0){$v$}}%
\put(30.9300,-7.1200){\makebox(0,0){$P$}}%
%
\special{sh 1.000}%
\special{ia 574 818 36 36  0.0000000 6.2831853}%
\special{pn 8}%
\special{ar 574 818 36 36  0.0000000 6.2831853}%
%
\special{sh 1.000}%
\special{ia 994 818 36 36  0.0000000 6.2831853}%
\special{pn 8}%
\special{ar 994 818 36 36  0.0000000 6.2831853}%
%
\special{sh 1.000}%
\special{ia 1554 818 36 36  0.0000000 6.2831853}%
\special{pn 8}%
\special{ar 1554 818 36 36  0.0000000 6.2831853}%
%
\special{pn 8}%
\special{pa 854 538}%
\special{pa 1414 1658}%
\special{fp}%
\special{pa 1694 538}%
\special{pa 1134 1658}%
\special{fp}%
%
\special{sh 1.000}%
\special{ia 1694 538 36 36  0.0000000 6.2831853}%
\special{pn 8}%
\special{ar 1694 538 36 36  0.0000000 6.2831853}%
%
\special{sh 1.000}%
\special{ia 854 538 36 36  0.0000000 6.2831853}%
\special{pn 8}%
\special{ar 854 538 36 36  0.0000000 6.2831853}%
%
\special{sh 1.000}%
\special{ia 1274 1378 36 36  0.0000000 6.2831853}%
\special{pn 8}%
\special{ar 1274 1378 36 36  0.0000000 6.2831853}%
%
\special{sh 1.000}%
\special{ia 1134 1658 36 36  0.0000000 6.2831853}%
\special{pn 8}%
\special{ar 1134 1658 36 36  0.0000000 6.2831853}%
%
\special{sh 1.000}%
\special{ia 1414 1658 36 36  0.0000000 6.2831853}%
\special{pn 8}%
\special{ar 1414 1658 36 36  0.0000000 6.2831853}%
%
\special{sh 1.000}%
\special{ia 1974 818 36 36  0.0000000 6.2831853}%
\special{pn 8}%
\special{ar 1974 818 36 36  0.0000000 6.2831853}%
%
\special{sh 1.000}%
\special{ia 2534 818 36 36  0.0000000 6.2831853}%
\special{pn 8}%
\special{ar 2534 818 36 36  0.0000000 6.2831853}%
%
\special{sh 1.000}%
\special{ia 2954 818 36 36  0.0000000 6.2831853}%
\special{pn 8}%
\special{ar 2954 818 36 36  0.0000000 6.2831853}%
%
\special{pn 8}%
\special{pa 2954 818}%
\special{pa 574 818}%
\special{fp}%
%
\special{pn 8}%
\special{pa 1448 880}%
\special{pa 574 880}%
\special{fp}%
%
\special{pn 8}%
\special{pa 1448 880}%
\special{pa 1064 1650}%
\special{fp}%
\special{sh 1}%
\special{pa 1064 1650}%
\special{pa 1112 1600}%
\special{pa 1088 1602}%
\special{pa 1076 1582}%
\special{pa 1064 1650}%
\special{fp}%
\put(10.9800,-13.7700){\makebox(0,0){$T$}}%
\put(19.7300,-9.5700){\makebox(0,0){$w^{(2)}_{2}$}}%
\put(25.3300,-9.5700){\makebox(0,0){$w^{(1)}_{2}$}}%
\end{picture}%
\caption{path $T$}
\label{fig4}
\end{center}
\end{figure}

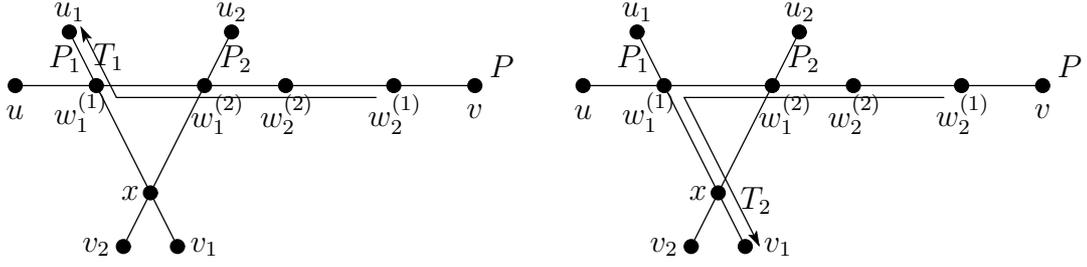
\begin{figure}
\begin{center}
\unitlength 0.1in
\begin{picture}( 54.6000, 13.2200)(  4.8300,-18.9200)
\put(35.1300,-11.6000){\makebox(0,0){$u$}}%
\put(38.5400,-11.5700){\makebox(0,0){$w^{(1)}_{1}$}}%
\put(45.6200,-11.7000){\makebox(0,0){$w^{(2)}_{1}$}}%
\put(46.3300,-6.3500){\makebox(0,0){$u_{2}$}}%
\put(37.9300,-6.3500){\makebox(0,0){$u_{1}$}}%
\put(39.3300,-18.6000){\makebox(0,0){$v_{2}$}}%
\put(45.2300,-18.6000){\makebox(0,0){$v_{1}$}}%
\put(41.0800,-15.7700){\makebox(0,0){$x$}}%
\put(46.5400,-8.7700){\makebox(0,0){$P_{2}$}}%
\put(37.7200,-8.7700){\makebox(0,0){$P_{1}$}}%
\put(58.9300,-11.5300){\makebox(0,0){$v$}}%
\put(60.3300,-9.1500){\makebox(0,0){$P$}}%
%
\special{sh 1.000}%
\special{ia 3514 1018 36 36  0.0000000 6.2831853}%
\special{pn 8}%
\special{ar 3514 1018 36 36  0.0000000 6.2831853}%
%
\special{sh 1.000}%
\special{ia 3934 1018 36 36  0.0000000 6.2831853}%
\special{pn 8}%
\special{ar 3934 1018 36 36  0.0000000 6.2831853}%
%
\special{sh 1.000}%
\special{ia 4494 1018 36 36  0.0000000 6.2831853}%
\special{pn 8}%
\special{ar 4494 1018 36 36  0.0000000 6.2831853}%
%
\special{pn 8}%
\special{pa 3794 738}%
\special{pa 4354 1858}%
\special{fp}%
\special{pa 4634 738}%
\special{pa 4074 1858}%
\special{fp}%
%
\special{sh 1.000}%
\special{ia 4634 738 36 36  0.0000000 6.2831853}%
\special{pn 8}%
\special{ar 4634 738 36 36  0.0000000 6.2831853}%
%
\special{sh 1.000}%
\special{ia 3794 738 36 36  0.0000000 6.2831853}%
\special{pn 8}%
\special{ar 3794 738 36 36  0.0000000 6.2831853}%
%
\special{sh 1.000}%
\special{ia 4214 1578 36 36  0.0000000 6.2831853}%
\special{pn 8}%
\special{ar 4214 1578 36 36  0.0000000 6.2831853}%
%
\special{sh 1.000}%
\special{ia 4074 1858 36 36  0.0000000 6.2831853}%
\special{pn 8}%
\special{ar 4074 1858 36 36  0.0000000 6.2831853}%
%
\special{sh 1.000}%
\special{ia 4354 1858 36 36  0.0000000 6.2831853}%
\special{pn 8}%
\special{ar 4354 1858 36 36  0.0000000 6.2831853}%
%
\special{sh 1.000}%
\special{ia 4914 1018 36 36  0.0000000 6.2831853}%
\special{pn 8}%
\special{ar 4914 1018 36 36  0.0000000 6.2831853}%
%
\special{sh 1.000}%
\special{ia 5474 1018 36 36  0.0000000 6.2831853}%
\special{pn 8}%
\special{ar 5474 1018 36 36  0.0000000 6.2831853}%
%
\special{sh 1.000}%
\special{ia 5894 1018 36 36  0.0000000 6.2831853}%
\special{pn 8}%
\special{ar 5894 1018 36 36  0.0000000 6.2831853}%
%
\special{pn 8}%
\special{pa 5894 1018}%
\special{pa 3514 1018}%
\special{fp}%
%
\special{pn 8}%
\special{pa 4038 1080}%
\special{pa 5382 1080}%
\special{fp}%
%
\special{pn 8}%
\special{pa 4038 1080}%
\special{pa 4424 1850}%
\special{fp}%
\special{sh 1}%
\special{pa 4424 1850}%
\special{pa 4412 1782}%
\special{pa 4400 1802}%
\special{pa 4376 1800}%
\special{pa 4424 1850}%
\special{fp}%
\put(49.2000,-11.7000){\makebox(0,0){$w^{(2)}_{2}$}}%
\put(54.8000,-11.6700){\makebox(0,0){$w^{(1)}_{2}$}}%
\put(44.1000,-16.1900){\makebox(0,0){$T_{2}$}}%
\put(5.7300,-11.6000){\makebox(0,0){$u$}}%
\put(9.1400,-11.5700){\makebox(0,0){$w^{(1)}_{1}$}}%
\put(16.2200,-11.7000){\makebox(0,0){$w^{(2)}_{1}$}}%
\put(16.9300,-6.3500){\makebox(0,0){$u_{2}$}}%
\put(8.5300,-6.3500){\makebox(0,0){$u_{1}$}}%
\put(9.9300,-18.6000){\makebox(0,0){$v_{2}$}}%
\put(15.5300,-18.6000){\makebox(0,0){$v_{1}$}}%
\put(11.6800,-15.7700){\makebox(0,0){$x$}}%
\put(17.1400,-8.7700){\makebox(0,0){$P_{2}$}}%
\put(8.3200,-8.7700){\makebox(0,0){$P_{1}$}}%
\put(29.5300,-11.5300){\makebox(0,0){$v$}}%
\put(30.9300,-9.1500){\makebox(0,0){$P$}}%
%
\special{sh 1.000}%
\special{ia 574 1018 36 36  0.0000000 6.2831853}%
\special{pn 8}%
\special{ar 574 1018 36 36  0.0000000 6.2831853}%
%
\special{sh 1.000}%
\special{ia 994 1018 36 36  0.0000000 6.2831853}%
\special{pn 8}%
\special{ar 994 1018 36 36  0.0000000 6.2831853}%
%
\special{sh 1.000}%
\special{ia 1554 1018 36 36  0.0000000 6.2831853}%
\special{pn 8}%
\special{ar 1554 1018 36 36  0.0000000 6.2831853}%
%
\special{pn 8}%
\special{pa 854 738}%
\special{pa 1414 1858}%
\special{fp}%
\special{pa 1694 738}%
\special{pa 1134 1858}%
\special{fp}%
%
\special{sh 1.000}%
\special{ia 1694 738 36 36  0.0000000 6.2831853}%
\special{pn 8}%
\special{ar 1694 738 36 36  0.0000000 6.2831853}%
%
\special{sh 1.000}%
\special{ia 854 738 36 36  0.0000000 6.2831853}%
\special{pn 8}%
\special{ar 854 738 36 36  0.0000000 6.2831853}%
%
\special{sh 1.000}%
\special{ia 1274 1578 36 36  0.0000000 6.2831853}%
\special{pn 8}%
\special{ar 1274 1578 36 36  0.0000000 6.2831853}%
%
\special{sh 1.000}%
\special{ia 1134 1858 36 36  0.0000000 6.2831853}%
\special{pn 8}%
\special{ar 1134 1858 36 36  0.0000000 6.2831853}%
%
\special{sh 1.000}%
\special{ia 1414 1858 36 36  0.0000000 6.2831853}%
\special{pn 8}%
\special{ar 1414 1858 36 36  0.0000000 6.2831853}%
%
\special{sh 1.000}%
\special{ia 1974 1018 36 36  0.0000000 6.2831853}%
\special{pn 8}%
\special{ar 1974 1018 36 36  0.0000000 6.2831853}%
%
\special{sh 1.000}%
\special{ia 2534 1018 36 36  0.0000000 6.2831853}%
\special{pn 8}%
\special{ar 2534 1018 36 36  0.0000000 6.2831853}%
%
\special{sh 1.000}%
\special{ia 2954 1018 36 36  0.0000000 6.2831853}%
\special{pn 8}%
\special{ar 2954 1018 36 36  0.0000000 6.2831853}%
%
\special{pn 8}%
\special{pa 2954 1018}%
\special{pa 574 1018}%
\special{fp}%
%
\special{pn 8}%
\special{pa 1098 1080}%
\special{pa 2442 1080}%
\special{fp}%
\put(19.8000,-11.7000){\makebox(0,0){$w^{(2)}_{2}$}}%
\put(25.4000,-11.6700){\makebox(0,0){$w^{(1)}_{2}$}}%
%
\special{pn 8}%
\special{pa 1098 1080}%
\special{pa 916 716}%
\special{fp}%
\special{sh 1}%
\special{pa 916 716}%
\special{pa 928 786}%
\special{pa 940 764}%
\special{pa 964 768}%
\special{pa 916 716}%
\special{fp}%
\put(10.5600,-8.7000){\makebox(0,0){$T_{1}$}}%
\end{picture}%
\caption{path $T_{i}$}
\label{fig5}
\end{center}
\end{figure}

Note that $|V(vPw^{(1)}_{2})|\leq |V(vPw^{(2)}_{2})|$.
Since the path $uP\check{w}^{(2)}_{1}$ contains no vertex in $V(P_{2})$, $T=uPw^{(2)}_{1}P_{2}v_{2}$ is a path in $G$ (see Figure~\ref{fig4}).
Since the length of $T$ is $(|V(uPw^{(1)}_{1})|-1)+(|V(Q_{1})|-1)+(|V(w^{(2)}_{1}P_{2}v_{2})|-1)$, $(|V(u_{2}P_{2}w^{(2)}_{1})|-1)+(|V(w^{(2)}_{1}P_{2}v_{2})|-1)=|V(P_{2})|-1=l(G)\geq (|V(uPw^{(1)}_{1})|-1)+(|V(Q_{1})|-1)+(|V(w^{(2)}_{1}P_{2}v_{2})|-1)$.
This together with (\ref{eq2.5}) leads to
\begin{align}
|V(u_{2}P_{2}w^{(2)}_{1})|\geq |V(uPw^{(1)}_{1})|+|V(Q_{1})|-1\geq |V(uPw^{(1)}_{1})|+f(G,\P).\label{eq2.10}
\end{align}
Since the path $\check{w}^{(1)}_{2}P\check{w}^{(1)}_{1}$ contains no vertex in $V(P_{1})$, both $T_{1}=\check{w}^{(1)}_{2}Pw^{(1)}_{1}P_{1}u_{1}$ and $T_{2}=\check{w}^{(1)}_{2}Pw^{(1)}_{1}P_{1}v_{1}$ are paths in $G$ (see Figure~\ref{fig5}).
Since the length of $T_{1}$ is $(|V(\check{w}^{(1)}_{2}Pw^{(1)}_{1})|-1)+(|V(w^{(1)}_{1}P_{1}u_{1})|-1)$, we have $(|V(vPw^{(1)}_{2})|-1)+(|V(w^{(1)}_{2}Pw^{(1)}_{1})|-1)+(|V(w^{(1)}_{1}Pu)|-1)=|V(P)|-1=l(G)\geq (|V(\check{w}^{(1)}_{2}Pw^{(1)}_{1})|-1)+(|V(w^{(1)}_{1}P_{1}u_{1})|-1)$.
Consequently, we have $|V(vPw^{(1)}_{2})|+|V(w^{(1)}_{1}Pu)|\geq |V(w^{(1)}_{1}P_{1}u_{1})|$.
By comparing the length of $P$ and $T_{2}$, we also have $|V(vPw^{(1)}_{2})|+|V(w^{(1)}_{1}Pu)|\geq |V(w^{(1)}_{1}P_{1}v_{1})|$.
Hence
\begin{eqnarray}
l(G) &=& |V(P_{1})|-1\nonumber \\
&=& |V(u_{1}P_{1}w^{(1)}_{1})|+|V(w^{(1)}_{1}P_{1}v_{1})|-2\nonumber \\
&\leq & 2(|V(vPw^{(1)}_{2})|+|V(w^{(1)}_{1}Pu)|)-2.\label{eq2.11}
\end{eqnarray}
Recall that the length of the unique $\{u\}$-$V(Q_{1})$ path on $P$ (i.e. $uPw^{(1)}_{1}$) is at least that of the unique $\{v\}$-$V(Q_{2})$ path on $P$ (i.e. $vPw^{(1)}_{2}$).
Hence $|V(uPw^{(1)}_{1})|\geq |V(vPw^{(1)}_{2})|$.
By (\ref{eq2.11}), $l(G)\leq 2(|V(vPw^{(1)}_{2})|+|V(w^{(1)}_{1}Pu)|)-2\leq 4|V(uPw^{(1)}_{1})|-2$, and so $|V(uPw^{(1)}_{1})|\geq (l(G)+2)/4$.
This together with (\ref{eq2.10}) implies that
\begin{align}
|V(u_{2}P_{2}w^{(2)}_{1})|\geq \frac{l(G)+2}{4}+f(G,\P).\label{eq2.12}
\end{align}
By (\ref{eq2.7}), (\ref{eq2.8}) and (\ref{eq2.12}),
\begin{eqnarray}
l(G) &=& |V(P_{2})|-1\nonumber \\
&=& |V(u_{2}P_{2}w^{(2)}_{1})|+|V(w^{(2)}_{1}P_{2}x)|+|V(xP_{2}v_{2})|-3\nonumber \\
&\geq & (\frac{l(G)+2}{4}+f(G,\P))+(2f(G,\P)+1)+2f(G,\P)-3\nonumber \\
&=& \frac{l(G)-6}{4}+5f(G,\P),\nonumber
\end{eqnarray}
and so
\begin{align}
l(G)\geq \frac{20f(G,\P)-6}{3}.\label{eq2.13}
\end{align}
By the choice of $P$, $t_{\P}(P')\geq 2$ for every $P'\in \P$.
By Lemma~\ref{lem2.2}, $\sum _{P'\in \P}|X_{\P}(P')|\geq \sum _{P'\in \P}t_{\P}(P')(f(G,\P)-1)\geq 6(f(G,\P)-1)$.
This together with Lemma~\ref{lem2.1} and (\ref{eq2.13}) implies that
\begin{eqnarray}
n &\geq & \frac{3l(G)+\sum _{P'\in \P}|X_{\P}(P')|+3}{2}\nonumber \\
&\geq & \frac{3\cdot \frac{20f(G,\P)-6}{3}+6(f(G,\P)-1)+3}{2}\nonumber \\
&=& \frac{26f(G,\P)-9}{2},\nonumber
\end{eqnarray}
and hence $f(G,\P)\leq (2n+9)/26$.

\medskip
\noindent
\textbf{Case 2:} $t_{\P}(P)\geq 3$.

By the choice of $P$, $t_{\P}(P')\geq 3$ for every $P'\in \P$.
By Lemma~\ref{lem2.2}, $\sum _{P'\in \P}|X_{\P}(P')|\geq \sum _{P'\in \P}t_{\P}(P')(f(G,\P)-1)\geq 9(f(G,\P)-1)$.
This together with Lemma~\ref{lem2.1} and (\ref{eq2.9}) implies that
\begin{eqnarray}
n &\geq & \frac{3l(G)+\sum _{P'\in \P}|X_{\P}(P')|+3}{2}\nonumber \\
&\geq & \frac{3(6f(G,\P)-2)+9(f(G,\P)-1)+3}{2}\nonumber \\
&=& \frac{27f(G,\P)-12}{2},\nonumber
\end{eqnarray}
and hence $f(G,\P)\leq (2n+12)/27$.

This completes the proof of Theorem~\ref{thm1}.
\qed

To conclude this section, we propose the following conjecture.

\begin{con}
\label{con2}
Let $G$ be a connected graph, and let $\P\subseteq \L(G)$ with $|\P|=3$.
If there exists a path $P\in \P$ with $t_{\P}(P)=2$, then $f(G,\P)=0$.
\end{con}

If Conjecture~\ref{con2} is true, then we can improve the upper bound of $f(G,\P)$ in Theorem~\ref{thm1} to $(2n+12)/27$ (by the argument in the proof of Theorem~\ref{thm1}).

\section{Bounding the value of $f(G,\P)$ by a sublinear function}

A function $g$ is \textit{sublinear}
if $\lim_{n \rightarrow +\infty} \frac{g(n)}{n} = 0$. 
It follows from the definition that, if $g$ is sublinear, then for any two constants $c_0, c_1$, we have $g(c_0t+c_1)<t$ for any large $t$. 
Here we pose the following new conjecture,
which concerns Conjecture \ref{Gallai}.
Although Conjecture \ref{sublinear} is 
seemingly weaker than Conjecture \ref{Gallai},
we will show that 
Conjecture~\ref{sublinear} is indeed equivalent with Conjecture~\ref{Gallai}.

\begin{con}
\label{sublinear}
There exists a sublinear non-decreasing function $g$ 
such that 
for every connected graph $G$ of order $n$
and every subset $\P$ of $\L(G)$ with $|\P|=3$,
$f(G, \P) \leq g(n)$.
\end{con}

To prove that this seemingly weaker conjecture is equivalent to 
Conjecture~\ref{Gallai}, we first show that for a given graph $G$ with a set 
$\{P_1,P_2, P_3\}$ of three longest paths one can choose a subdivision of $G$ so that 
subdivisions of $P_i$'s $i=1,2,3$ are the new longest paths and show that the 
minimum distance from these three subdivided paths in the subdivided graph 
grows linearly in the order of subdivision. For the exact statement we introduce 
the following notation.

Let $G$ be a connected graph and let $\mathcal P=\{P_1,P_2, P_3\}$  be a set of three 
longest paths. Let $G'$ be obtained by adding a new edge to each end-vertex 
of  $P_i$'s, $i=1,2,3$; thus, minimum of two and maximum of six new vertices and 
edges are added. Let $P'_i$, $i=1,2,3$ be the path corresponding to $P_i$ with two 
new edges at the two ends.  We define $G^t$ to be the graph obtained from $G'$ 
by subdividing 
each edge $t$ times. Let $P^t_i$, $i=1,2,3$ be the corresponding path of $P'_i$ 
in $G^t$.  We write 
$\mathcal P^{t}= \{P^t_1,P^t_2, P^t_3\}$. Also, let $V_{f(G,\P)}=\{v\in V(G)\mid \sum_{P\in \P}d_{G}(v,V(P))=f(G,\P)\}$.

We have the following proposition.

\begin{prop}
 Given a connected graph $G$ and a set $\mathcal P=\{P_1,P_2, P_3\}$ of three longest 
paths, the set $\mathcal P^{t}= \{P^t_1,P^t_2, P^t_3\}$ is a set of three 
longest paths of $G^t$. Furthermore, $f(G^t, \P^t)=(t+1)f(G, \P)$.
\end{prop}

\proof
 The first assertion is easy to check. To prove the second assertion, we show that a 
vertex  of $V_{f(G^t,\mathcal{P}^t)}$ could be chosen as 
an original vertex of $G$. The assertion then would follow, as the vertex of $G$ attaining  
the distance sum $f(G, \P)$ of $\P$ satisfies $(t+1)f(G, \P)$ for the distance sum of $\P^t$ in $G^t$ as well.

Now let $u$ be a vertex attaining the distance sum $f(G^t, \P^t)$ from $\P^t$. 
It is easy to check that $u$ is not an end-vertex of $P^t_i$ for any $i$. If $u\in 
V(G)$, then we have nothing to prove. Otherwise $u$ is a new vertex 
subdividing an edge, say $xy$, of $G$. If all the shortest paths from $u$ to 
$P^t_i$, $i=1,2,3$, go through $x$ (or $y$) then replacing $u$ by $x$ (or $y$) provides a 
smaller distance sum than $f(G^t, \P^t)$, a contradiction. Thus we may assume, without loss of 
generality, that two of the shortest paths from $u$ to $P^t_i$ go through $x$ 
and the third one goes through $y$. 
In such a case again by replacing $u$ by $x$ we 
will 
have a smaller distance sum than $f(G^t, \P^t)$, a contradiction. 
We note that if $u$ belongs to one or two
of these paths then so are $x$ and $y$, thus this would not affect the 
argument. The contradiction proves that $u$ must 
be a vertex of $G$ and we have $f(G^t, \P^t)=(t+1)f(G, \P)$.
\qed

Keeping the above proposition in mind, we can prove the following theorem. 

\begin{thm}\label{order}
Conjecture~\ref{Gallai} is true if and only if 
Conjecture~\ref{sublinear} is true.
\end{thm}

\proof
The ``only if'' part is trivial,
and hence we only show the ``if'' part.

Suppose that $G$ together with $\P =\{P_1,P_2, P_3\}$ is a counterexample for 
Conjecture~\ref{Gallai}, i.e.,  $f(G, \P)\geq 1$. The subgraph of $G$ 
induced by edges and vertices of $P_1,P_2,P_3$  is also a
counterexample (where $\P$ is also a set of non-intersecting three longest 
paths). Note that, in view of Proposition~\ref{prop1}, such a subgraph is connected.  
Thus we may assume from the start that vertices and edges of $G$ are union of 
vertices and edges of $P_1,P_2,P_3$. Let $n_0$ be the number of vertices of $G$. 
Since $G$ is union of three paths each of length at most $n_0-1$, we conclude 
that $G$ has at most $3(n_0-1) (<3(n_0+1))$ edges. 

Hence, by the construction of $G^t$, we have $|V(G^t)|\leq n_0+3(n_0+1)t+6$. 
On the other hand, we have 
$f(G^t, \P^t)=(t+1)f(G, \P)\geq t$.  Hence for constants $c_0=3n_0+3$ and $c_1=n_0+6$
we have $g(c_0t+c_1)\geq t$ (because $g$ is non-decreasing). 
However, this contradicts the fact that $g$ is a sublinear function. 
\qed

In conclusion, Theorem \ref{order} tells us that giving a substantial improvement on the magnitude of the upper bound of $f(G,\P)$ in Theorem~\ref{thm1} settles the longstanding conjecture on intersecting three longest paths in a connected  graph. \\

\noindent\textbf{Acknowledgments}

We would like to thank Dr.Valentin Borozan for a fruitful discussion concerning this paper. 
The first author would like to thank the laboratory LRI of the University Paris South and Digiteo foundation for their generous hospitality.  He was able to carry out part of this research during his visit there. Also, the first author's research is supported by the Japan Society for the Promotion of Science Grant-in-Aid for Young Scientists (B) (20740095).  
The second author's research is in part supported by
the Japan Society for the Promotion of Science Grant-in-Aid for Young Scientists (B) (26800086).  
The fourth author's research is in part supported by
the Japan Society for the Promotion of Science Grant-in-Aid for Young Scientists (B) (25871053),
and by Grant for Basic Science Research Projects from The Sumitomo Foundation.

\end{document}